\documentclass{amsart}
\usepackage{fullpage,graphicx,psfrag,epsfig,amsmath,amsfonts,amssymb,fancybox,amsthm}
\usepackage[dvips]{color}
\newcommand{\eq}{\begin{equation}}
\newcommand{\en}{\end{equation}}
\newcommand{\eqa}{\begin{eqnarray}}
\newcommand{\ena}{\end{eqnarray}}
\newcommand{\eqas}{\begin{eqnarray*}}
\newcommand{\enas}{\end{eqnarray*}}

\definecolor{Red}{rgb}{1,0,0}
\definecolor{Blue}{rgb}{0,0,1}
\definecolor{Green}{rgb}{0,1,0}
\definecolor{Yellow}{rgb}{1,1,0}
\definecolor{Cyan}{rgb}{0,1,1}
\definecolor{Magenta}{rgb}{1,0,1}
\definecolor{Orange}{rgb}{1,.5,0}
\definecolor{Violet}{rgb}{.5,0,.5}
\definecolor{Purple}{rgb}{.75,0,.25}
\definecolor{Brown}{rgb}{.75,.5,.25}
\definecolor{Grey}{rgb}{.5,.5,.5}

\newtheorem{theorem}{Theorem}[section]
\newtheorem{lemma}[theorem]{Lemma}

\newtheorem{coro}[theorem]{Corollary}
\theoremstyle{definition}
\newtheorem{definition}[theorem]{Definition}
\newtheorem{example}[theorem]{Example}

\theoremstyle{remark}
\newtheorem{rem}[theorem]{Remark}

\numberwithin{equation}{section}

\begin{document}
\title{The correlation decay (CD) tree and strong spatial mixing in multi-spin
systems}
\author{Chandra Nair}
\address{Microsoft Research, Redmond, WA 98052 }
\email{cnair@microsoft.com}
\author{Prasad Tetali}
\address{School of Mathematics and college of computing, Georgia Institute of Technology, Atlanta, GA 30332}
\email{tetali@math.gatech.edu}
\thanks{Part of this work was done while Prasad Tetali was visiting Microsoft Research during 2006.}


\begin{abstract}
This paper deals with the construction of a {\em correlation
decay} tree (hypertree) for interacting systems modeled using
graphs (hypergraphs) that can be used to compute the marginal
probability of any vertex of interest. Local message passing
equations have been used for some time to approximate the marginal
probabilities in graphs but it is known that these equations are
incorrect for graphs with loops. In this paper we construct, for
any finite graph and a fixed vertex, a finite tree with
appropriately defined boundary conditions so that the marginal
probability on the tree at the vertex matches that on the graph.
For several interacting systems, we show using our approach that
if there is very strong spatial mixing on an infinite regular
tree, then one has strong spatial mixing for any given graph with
maximum degree bounded by that of the regular tree. Thus we
identify the regular tree as the worst case graph, in a weak
sense, for the notion of strong spatial mixing.
\end{abstract}

\maketitle

\section{Introduction}

In this paper we show that computation of the marginal probability for
a vertex in a graphical model can be reduced to the computation of
the marginal probability of the vertex in a rooted tree of self-avoiding
walks, with appropriately defined boundary conditions. The
computation tree approach for graphical models has been used by
\cite{wei06}, \cite{bag06}, \cite{gak07}, for the problems of
independent sets, colorings and list-colorings. In \cite{jus06}, the
work of \cite{wei06} for computing marginal probabilities was
extended to inference problems in general two spins models.
 Our work builds on \cite{wei06, jus06}, and
demonstrates how the computation tree can be extended to more than
two spins and also for more than two-body interactions. This leads
to a different tree (the {\em
correlation decay} tree), which in a sense is more natural than the dynamic programming based tree
of \cite{gak07} for the case of multiple spins. Further, this
approach also yields a tree for the case of multi-spin
interactions with multiple spins.

A practical motivation for the creation of a tree structure is the
following. The feasible algorithms for computation of marginal
probabilities in large interacting systems are constrained to be
distributed and local. This requirement has given rise to message
passing algorithms (like belief propagation) for systems modeled
using graphs. Unfortunately, these algorithms do not necessarily
give the correct answer for graphs with many loops, and may not
even converge. However, for a tree it is known that the equations
are exact and the marginal probability at the root can be computed
in a single iteration by starting from the leaves. Thus, if for
any graph one can show the existence of a tree, that respects the
locality, in which the same marginal probability results, then one
can use the exactness of the message passing algorithms on a tree
to obtain a convergent, distributed, local algorithm for the
computation of marginal probabilities on the original graph.

The caveat with this approach is that the size of the tree can be
exponentially large compared to the original graph. So even though
the computations are exact, they may not be efficient in practice.
However, for certain interesting counting problems
\cite{wei06,gak07,bgknt06} approximation algorithms have been
designed using the notion of spatial correlation decay, where the
influence of the boundary at a root decays as the spatial distance
between the boundary and the root increases.
Hence pruning the tree to an efficiently computable
neighborhood usually yields good and efficient approximations.
Thus, to design efficient algorithms it would be useful to show
some kind of decay of correlation in the tree structure that is
presented here (and hence the name {\em correlation decay} tree).

The second part of this paper addresses this issue of spatial
correlation decay. We show that, for lots of systems of interest,
if there is ``very strong spatial mixing" in the infinite regular
tree of degree $D$, then there also exists ``strong spatial
mixing" for any graph with maximum degree $D$. So, in a loose
sense, the infinite regular tree is indeed a worst case graph for
correlation decay. The fact that some form of strong spatial
mixing in the infinite regular tree should imply strong spatial
mixing in graphs for a general multi-spin system was  conjectured
by E. Mossel, \cite{mos07}. (In the case of independent sets and
colorings, the infinite tree being the worst case for the onset of
multiple Gibbs measures was conjectured by A. Sokal \cite{sok00}.)

In the next section, we prove the generalization of the result in
\cite{wei06} to the case of multiple-spins but still restricting
ourselves to two-body (pairwise) interactions.

\section{Preliminaries}

Consider a finite spin system with pairwise interactions, and
modeled as a graph, $G=(V,E)$. Let the partition function of this
spin system be denoted by
\[ Z_G = \sum_{\vec{x} \in X^n}
\prod_{(i,j) \in G} \Phi_{i,j}(x_i,x_j) \prod_{i \in V}
\phi_i(x_i). \] Let $\Lambda \subseteq [n]$ be a subset of {\em
frozen} vertices (i.e. vertices whose spin values are fixed) and
let
\[ Z_G^\Lambda =
\sum_{\vec{x} \notin X_\Lambda} \prod_{(i,j) \in
G}\Phi_{i,j}(x_i,x_j) \prod_{i \in V} \phi_i(x_i).  \] We wish to
compute the following marginal probability with respect to the
Gibbs measure,
\begin{equation}
\label{eq:GibbsMargProb1} P_G(x_1 = \sigma | X_\Lambda ) =
\frac{1}{Z_G^\Lambda} \sum_{\stackrel{x_1 = \sigma,} {\vec{x}
\notin X_\Lambda}} \prod_{(i,j) \in G}\Phi_{i,j}(x_i,x_j) \prod_{i
\in V} \phi_i(x_i).
\end{equation}
Instead of performing this marginal probability computation in the
original graph $G$ we shall create a {\em correlation decay} (CD)
tree, $T_\Lambda$, on which the same marginal probability results
by performing the computation as described in Section
\ref{sse:compsec}.

\subsection{The CD Tree}

\label{sse:comptree}

Similar tree constructions can be found for restricted classes of
spin systems by \cite{wei06,fis59,gmp04,scs05,bag06,jus06}, and in
particular the one in \cite{wei06}. Our starting point of the
tree is  the same as in \cite{wei06}, i.e. we begin by labeling
the edges of the graph; draw the tree of self avoiding walks,
$T_{saw}$; and include the vertices that close a cycle. In
\cite{wei06}, the vertices that close a cycle were denoted as
occupied or unoccupied depending on whether the edge closing the
cycle in $T_{saw}$ was larger than the edge beginning the cycle or
not.

Our main point of deviation from the construction in \cite{wei06}
is in the treatment of vertices that close the cycle that were
appended in $T_{saw}$. The vertices that close the cycle with
higher numbered edges than those that begin the cycle (i.e. those
that were marked occupied) are now constrained to take a
particular spin value $\sigma_q$.    The vertices that close the
cycle with lower numbered edges (i.e. the unoccupied vertices) are
constrained to take the same value as the occurrence of it earlier
in the graph, i.e. the value of the vertex that begins the cycle.
This constraint is denoted by a {\em coupling line} and influences
the way the marginal probabilities are computed on the tree. The
tree thus obtained is called the CD-tree, $T_{CD}$, associated
with graph $G$.

\begin{definition}
\label{def:coupling} A {\em coupling line} on a rooted tree is a
virtual line connecting a vertex $u$ to some vertex $v$ in the
subtree below $u$. This line will play a role in the computation of the
marginal probabilities as will be explained in detail later. In brief words,
when one descends into the subtree of $u$ to compute the marginal probability
that $u$ assumes a spin $\sigma_i$, then the vertex $v$ becomes frozen to
$\sigma_i$, the same as $u$. Thus, the spin to which $v$ is frozen is coupled to
the spin of $u$, whose marginal probability is being determined.
\end{definition}

\begin{rem}
\label{rem:topofcl} One can easily make the following observations
regarding coupling lines. A vertex can be the top end point of
several coupling lines and indeed the number of coupling lines
from any point is related to the number of cycles the vertex is
part of in a certain subgraph of the original graph. A vertex can
only be the bottom end point of a unique coupling line and for
every such point, there is a unique twin point whose spin is
frozen to $\sigma_q$, corresponding to traversing the cycle in the
opposite direction.
\end{rem}

\subsection{Computation of marginal probabilities on the CD
tree} \label{sse:compsec} Here we describe the algorithm for
computing the marginal probability at the root for a tree with
coupling lines. Let $T$ be a rooted tree with frozen vertices $\Lambda$. In
the tree presented in the previous section, the set $\Lambda$ is
also assumed to contain the vertices frozen to $\sigma_q$. Consider
the recursion
\begin{equation}
\label{eq:treerec0} R_T^{\sigma_\Lambda}(\sigma_v) =
\frac{\phi_v(\sigma_v)}{\phi_v(\sigma_q)} \prod_{i=1}^d
\frac{\sum_{l=1}^q \Phi_{v,u_i}(\sigma_v,\sigma_l)
R_{T_i}^{\sigma_{\Lambda_i}}(u_i=\sigma_l)}{\sum_{l=1}^q
\Phi_{v,u_i}(\sigma_q,\sigma_l)
R_{T_i}^{\sigma_{\Lambda_i}}(u_i=\sigma_l)}\,.
\end{equation}
At this step (proven by the next theorem) we will be computing the
ratio of the probability that the root assumes a spin $\sigma_v$ (with respect
to the reference spin
$\sigma_q$), and therefore the lower end points of the coupling lines joined to
the root to be frozen to
$\sigma_v$. Thus the set of frozen vertices $\Lambda$ gets appended
with this subset of vertices; and the subset of this enhanced
$\Lambda$ that is in the subtree of the $i$th child is denoted as
$\sigma_{\Lambda_i}$. (There is an abuse of notation in that
$\sigma_{\Lambda_i}$ depends on the spin $\sigma_v$ as $\Lambda$
gets appended with the new vertices frozen by the dotted lines to
$\sigma_v$.) One can use the above recursion to recursively compute
the ratios for the correlation decay tree. The validity of this computation
forms the basis of the next theorem.

\begin{rem}
Consider a rooted tree with $D$ denoting the maximum number of
children for any vertex. Let $C$ denote the computation time
required for one step of the recursion in \eqref{eq:treerec0}, then
it is clear that computing the probability at the root given the
marginal probabilities at depth $\ell$ requires
$\Theta([(q-1)D]^\ell)$ time. The hidden constants in $\Theta$
depend on $C$ and $q$. Observe that a bound for the computation
time, $t_\ell$, at depth $\ell$ can be obtained via the recursion
$t_\ell \leq qC + [(q-1) D] t_{\ell - 1}$.
\end{rem}

Note that whenever the tree visits a frozen vertex, the subtree
under the frozen vertex can be pruned as this does not affect the
computation. Similarly the subtree under a vertex that is also
below the lower end of the virtual coupling line can be pruned.
This leads to a subtree, $T_{CD}^\Lambda$, of $T_{CD}$.

\begin{example}
We shall demonstrate this construction and computation using the
following example graph with edges labeled in the usual
lexicographic order. We shall retain the labeling of vertices on
$T_{CD}$ to reflect its origin from $G$ but other than that they
play no role in spin assignments and two similarly labeled
vertices can have arbitrary spin assingments in general.

\begin{figure}[ht]
\begin{center}
\includegraphics[width=0.5\linewidth,angle=0]{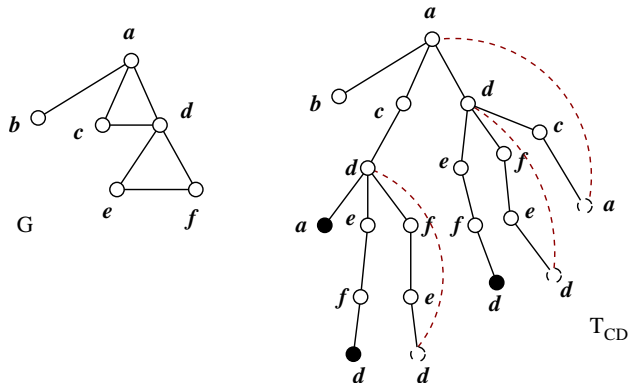}
\caption{The construction of the CD tree: The light dotted lines in
the figure denote the virtual {\em coupling lines}.}
\label{fig:ct}
\end{center}
\end{figure}

Let us assume that we are interested in computing the marginal
probability of the vertex $a$ for valid 5-colorings of the graph
$G$ using the tree $T_{CD}$ on the right. A coloring is valid if
no two adjacent vertices are assigned the same color. For this
interacting system $\sigma \in \{1,2,3,4,5\}$ and
\begin{equation*}
 \Phi(\sigma_i,\sigma_j) = \left\{ \begin{array}{ll} 1 &
\mbox{if}~ \sigma_i \neq \sigma_j~\\ 0 & \mbox{if}~ \sigma_i = \sigma_j
\end{array} \right.
\end{equation*}
 and the potential function $\phi(\sigma_i) = 1$. Let us assume that the
vertices $b,c$ are frozen to spins $2,3$ respectively and the
reference spin $\sigma_q = 4$. It is easy to see using symmetry or
explicit computation that $a$ takes spins $1,4,5$ with probability
$1/3$ each, or in other words the ratios (with respect to color
4), $R_G(1) = R_G(4) = R_G(5) = 1$. The pruned subtree
$T_{CD}^\Lambda$ can be drawn as in Figure \ref{fig:pt}.
\begin{figure}[ht]
\begin{center}
\includegraphics[width=0.45\linewidth,angle=0]{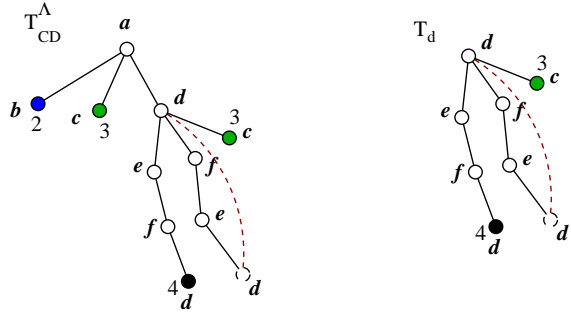}
\caption{The CD-tree $T_{CD}^\Lambda$ and the subtree $T_d$,
pruned by the frozen vertices and coupling lines. The frozen
colors are written adjacent to vertices in $\Lambda$. }
\label{fig:pt}
\end{center}
\end{figure}

Equation \eqref{eq:treerec0} gives
\begin{equation}
\label{eq:comp11} R_{T_{CD}}(1) = \frac{R_{T_d}(2) + R_{T_d}(4) +
R_{T_d}(5)}{R_{T_d}(1) + R_{T_d}(2) + R_{T_d}(5)},
\end{equation}
where $T_d$ represents the subtree of $T_{CD}^\Lambda$ under
vertex $d$. The frozen subtrees $T_d$ for the four computations
$R_{T_d}(1)$, $R_{T_d}(2)$, $ R_{T_d}(4), R_{T_d}(5)$ are
represented in Figure \ref{fig:pt2}.

\begin{figure}[ht]
\begin{center}
\includegraphics[width=0.55\linewidth,angle=0]{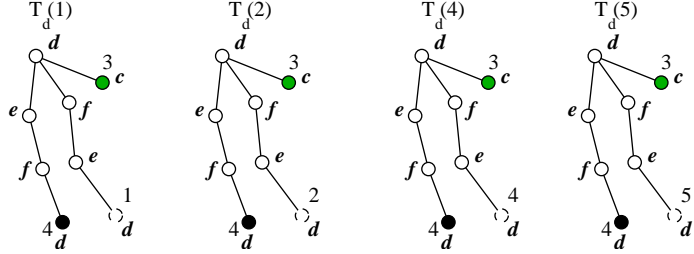}
\caption{The subtree $T_d(\cdot)$ for the computations
$R_{T_d}(1), R_{T_d}(4), R_{T_d}(4), R_{T_d}(5),$ respectively.
Note the spin of the new frozen vertex as forced by the coupling
line in the four cases.} \label{fig:pt2}
\end{center}
\end{figure}

The resultant subtrees $T_d(\cdot)$ have the usual computation
procedure (i.e. they do not have coupling lines);  for example, the
value $R_{T_d}(1)$ can be computed as
\begin{equation*}
\begin{split}
 R_{T_d}(1) & = \left(\frac{R_{T_e}(2) + R_{T_e}(3) + R_{T_e}(4)
+ R_{T_e}(5)}{R_{T_e}(1) + R_{T_e}(2) + R_{T_e}(3) +
R_{T_e}(5)}\right)\left( \frac{R_{T_f}(2) + R_{T_f}(3) +
R_{T_f}(4) +
R_{T_f}(5)}{R_{T_f}(1) + R_{T_f}(2) + R_{T_f}(3) + R_{T_f}(5)}\right) \\
& = \left(\frac{ \frac34 + \frac34 + \frac34 + \frac34}{1 +
\frac34 + \frac34 + \frac34}\right)\left( \frac{ \frac34 + \frac34
+ 1 + \frac34}{\frac34 + \frac34 + \frac34 + \frac34}\right)  = 1.
\end{split}
\end{equation*}
By symmetry to the previous computation $R_{T_d}(2) =
R_{T_d}(5)=1$ and from the definition, $R_{T_d}(4)=1$. Thus from
\eqref{eq:comp11} one obtains
\[ R_{T_{CD}}(1) = \frac{ 1 + 1 + 1}{ 1 + 1 + 1} = 1,\]
as desired.

\end{example}

\begin{rem}
The next theorem and its proof is essentially the same as in
\cite{wei06a}; therefore we will use the same notation whenever
possible and skip the details of similar arguments.
\end{rem}

\begin{theorem}
\label{th:comptreepresmargprob} For every graph $G=(V,E)$, every
$\Lambda \subseteq V$, any configuration $ \sigma_\Lambda$, and
all $\sigma_v$
\[ R_G^{\sigma_\Lambda}(v=\sigma_v) =
\mathbb{R}_{T_{CD}}^{\sigma_\Lambda}(v=\sigma_v), \] where
$\mathbb{R}_{T_{CD}}^{\sigma_\Lambda}(v=\sigma_v)$ stands for the
ratio (with respect to the reference spin, say $q$) of the
probability that the root $v$ of $T_{CD}$ has spin $\sigma_v$ when
the computation is performed as described above. The actual
probabilities can be computed from the ratios by normalizing them
such that the probabilities sum to one.
\end{theorem}

\begin{proof}
Let $\sigma_q$ be a fixed spin. Define the ratios
$$R_G^{\sigma_\Lambda}(\sigma_v)
\stackrel{\triangle}{=}
\frac{p_G^{\sigma_\Lambda}(v=\sigma_v)}{p_G^{ \sigma_\Lambda }
(v=\sigma_q) }. $$ Let $d$ be the degree of vertex $v$ and let
$u_i, 1 \leq i \leq d$ be its neighbors. If the graph $G$ was
indeed a tree $T$, then we can see that the following exact
recursion
\begin{equation}
\label{eq:treerec} R_T^{\sigma_\Lambda}(\sigma_v) =
\frac{\phi_v(\sigma_v)}{\phi_v(\sigma_q)} \prod_{i=1}^d
\frac{\sum_{l=1}^q \Phi_{v,u_i}(\sigma_v,\sigma_l)
R_{T_i}^{\sigma_{\Lambda_i}}(u_i=\sigma_l)}{\sum_{l=1}^q
\Phi_{v,u_i}(\sigma_q,\sigma_l)
R_{T_i}^{\sigma_{\Lambda_i}}(u_i=\sigma_l)},
\end{equation} would hold, where $T_i$ is the
subtree associated with the neighbor $u_i$ obtained by removing
the $i$th edge of $v$, and $\sigma_{\Lambda_i}$ is the restriction
of $\sigma_\Lambda$ to $\Lambda \cap T_i$ and appended with the
new vertices frozen to $\sigma_v$ corresponding to the lower
endpoints of coupling lines originating from $v$.

Fixing the vertex of interest $v$, define $G'$ as the graph
obtained by making $d$ copies of the vertex $v$ and each $v_i$
having a single edge to $u_i$. In addition, the vertex potential
$\phi_v(\sigma_v)$ is re-defined to $\phi_v^{1/d}(\sigma_v)$. It
is easy to see that the following two ratios are equal
$$ \frac{p_G^{\sigma_\Lambda}(v=\sigma_v)}{p_G^{ \sigma_\Lambda
} (v=\sigma_q)} =
\frac{p_{G'}^{\sigma_\Lambda}(v_1=\sigma_v,...,v_d =
\sigma_v)}{p_{G'}^{ \sigma_\Lambda} (v_1=\sigma_q, ... , v_d =
\sigma_q)}.$$ Defining $$R_{G',v_i}^{\sigma_\Lambda
\tau_i}(\sigma_v) =
\frac{p_{G'}^{\sigma_\Lambda}(v_1=\sigma_v,...,v_i = \sigma_v,
v_{i+1} = \sigma_q, .. , v_d = \sigma_q)}{p_{G'}^{ \sigma_\Lambda}
(v_1=\sigma_v,...,v_{i-1} = \sigma_v, v_i = \sigma_q, .. , v_d =
\sigma_q)}$$
one sees that
$$ R_G^{\sigma_\Lambda}(\sigma_v) = \prod_{i=1}^d R_{G',v_i}^{\sigma_\Lambda
\tau_i}(\sigma_v).$$
It is easy to see that
$R_{G',v_i}^{\sigma_\Lambda \tau_i}(\sigma_v)$ is the ratio of the
probaility that the vertex $v_i = \sigma_v$ to the probability of
$v_i = \sigma_q$, conditioned on $\sigma_\Lambda$ and $\tau_i$,
 where $\tau_i$ denotes the configuration where vertices $v_1,...,v_{i-1}$ are
frozen to $\sigma_v$ and vertices $v_{i+1},...,v_{d}$ are frozen
to $\sigma_q$.

In $G'$, the vertex $v_i$ is only connected to $u_i$; and let
$G'\setminus v_i$ denote the connected component of $G'$ that
contains $u_i$ after the removal of the edge $(v_i,u_i)$.
Therefore
$$ R_{G',v_i}^{\sigma_\Lambda \tau_i}(\sigma_v) =
\frac{\phi_v^{1/d}(\sigma_v)}{\phi_v^{1/d}(\sigma_q)}
\frac{\sum_{l=1}^q \Phi_{v,u_i}(\sigma_v,\sigma_l) R_{G'\setminus
v_i}^{\sigma_{\Lambda} \tau_i}(u_i=\sigma_l)}{\sum_{l=1}^q
\Phi_{v,u_i}(\sigma_q,\sigma_l) R_{G'\setminus
v_i}^{\sigma_{\Lambda} \tau_i}(u_i=\sigma_l)}, $$ and hence
\begin{equation}
\label{eq:graphrecur} R_G^{\sigma_\Lambda}(\sigma_v) =
\frac{\phi_v(\sigma_v)}{\phi_v(\sigma_q)} \prod_{i=1}^d
\frac{\sum_{l=1}^q \Phi_{v,u_i}(\sigma_v,\sigma_l) R_{G'\setminus
v_i}^{\sigma_{\Lambda} \tau_i}(u_i=\sigma_l)}{\sum_{l=1}^q
\Phi_{v,u_i}(\sigma_q,\sigma_l) R_{G'\setminus
v_i}^{\sigma_{\Lambda} \tau_i}(u_i=\sigma_l)}.
\end{equation}
Observe that the recursion \eqref{eq:graphrecur} terminates since
at each step the number of unfixed vertices reduces by one.

\begin{rem}
Observe that the equation in \eqref{eq:graphrecur} is similar to
the one for the tree \eqref{eq:treerec}. This similarity will help
us identify the recursion \eqref{eq:graphrecur} to be exactly the
same one in $T_{CD}$ with the condition corresponding to
$\sigma_\Lambda$ along with the coupling of the values of vertices
that was used in its definition. The key difference between the
binary spin model in \cite{wei06} and this proof also lies here;
that in the binary spin model one of the spins was always the
reference spin and the other was the subject of the recursion.
Thus the coupling of the spin to its parent in $T_{CD}$ was
implicit.
\end{rem}

From the similarity of \eqref{eq:graphrecur} and
\eqref{eq:treerec}, one can use induction to complete the proof
provided that the graph $G'\setminus v_i$  with the condition
corresponding to $\sigma_\Lambda \tau_i$ leads to the same subtree
of $T_{CD}$ corresponding to the $i$-th child of the original root
with the condition corresponding to $\sigma_{\Lambda_i}$. It is
easy to observe that the two trees are the same -- both are paths
in $G$ starting at $u_i$, and copies of $v$ are  set to $\sigma_v$
if it is reached via a smaller numbered edge and set to $\sigma_q$
else. The above observation along with the fact that the stopping
rules coincide for the two recursions completes the proof of
Theorem \ref{th:comptreepresmargprob} using induction.

\end{proof}

\subsection{Multi-spin interactions}

\label{sse:mulspiter}

In this section, we extend the results of the previous section
from pairwise interactions to multi-spin interactions. The
underlying model can be depicted by a hypergraph with the
hyperedges denoting the vertices involved in an interaction.

Consider a finite spin system whose interactions can be modeled
as a hypergraph, $G=(V,E)$. Let the partition function of this
spin system be denoted by
\[ Z_G = \sum_{\vec{x} \in X^n}
\prod_{e \in E} \Phi_{e}(\vec{x}_e) \prod_{i \in V} \phi_i(x_i).
\]
As before, let $\Lambda \subseteq [n]$ be a subset of {\em frozen}
vertices (i.e. vertices whose spin values are fixed) and let
\[ Z_G^\Lambda =
\sum_{\vec{x} \notin X_\Lambda} \prod_{e \in E}\Phi_{e
}(\vec{x}_e) \prod_{i \in V} \phi_i(x_i).  \] We wish to compute
the following marginal probability with respect to the Gibbs
measure,
\begin{equation}
\label{eq:GibbsMargProb} P_G(x_1 = \sigma | X_\Lambda ) =
\frac{1}{Z_G^\Lambda} \sum_{\stackrel{x_1 = \sigma,} {\vec{x}
\notin X_\Lambda}} \prod_{e \in E}\Phi_{e}(\vec{x}_e) \prod_{i \in
V} \phi_i(x_i).
\end{equation}

\subsection{CD hypertrees on hypergraphs}

\label{sse:comptreehg} The motivation for the following hypertree
essentially comes from the proof of the CD tree in the previous
section. Let the $n$ vertices in  $V$ be numbered in some fixed order,
${1,...,n}$. The tree is constructed in a top down approach
just as the tree of self avoiding walks.

The procedure described below is similar to a generalization of
the tree of self avoiding walks for graphs. For ease of exposition
we will keep describe the construction using the following
example. Let $V = \{1,2,3,4,5\}$ and let the hyperedges be
$\{(1,2,3),(1,2,5),(1,3,4),(2,5,4)\}.$ Let us assume that vertex
$1$ is the root. From $G$ construct the graph $G_1$ with vertex
$1$ replicated thrice (equal to its degree) to $1_a,1_b,1_c$. Let
the resulting hyperedges be
$\{(1_a,2,3),(1_b,2,5),(1_c,3,4),(2,5,4)\}$. Observe that,
\begin{align*}
\frac{P_G(x_1 = \sigma_1)}{P_G(x_1= \sigma_0)}  = ~&
\frac{P_{G_1}(x_{1_a} = \sigma_1,
x_{1_b}=\sigma_1,x_{1_c}=\sigma_1)}{P_{G_1}(x_{1_a} = \sigma_0,
x_{1_b}=\sigma_0,x_{1_c}=\sigma_0)} \\
= ~& \frac{P_{G_1}(x_{1_a} = \sigma_1|
x_{1_b}=\sigma_0,x_{1_c}=\sigma_0)}{P_{G_1}(x_{1_a} = \sigma_0|
x_{1_b}=\sigma_0,x_{1_c}=\sigma_0)} \times \frac{P_{G_1}(x_{1_b} =
\sigma_1| x_{1_a}=\sigma_1,x_{1_c}=\sigma_0)}{P_{G_1}(x_{1_b} =
\sigma_0|
x_{1_a}=\sigma_1,x_{1_c}=\sigma_0)}  \\
& \quad \times \frac{P_{G_1}(x_{1_c} = \sigma_1|
x_{1_a}=\sigma_1,x_{1_b}=\sigma_1)}{P_{G_1}(x_{1_c} = \sigma_0|
x_{1_a}=\sigma_1,x_{1_b}=\sigma_1)}.
\end{align*}

Now consider a graph $H$ where vertex 1 has degree three and such
that the removal of vertex 1 and the three hyperedges, disconnects
the graph into 3 disconnected components. The first component,
$H_1$, contains the set of vertices $\{2^{(1)},3^{(1)},4^{(1)},5^{(1)},
1_b^{(1)}, 1_c^{(1)}\}$, with the vertices $1_b^{(1)}$ and
$1_c^{(1)}$ frozen to have spin $\sigma_0$. The hyperedges that
form part of this component (along with the root) are
$\{(1,2^{(1)},3^{(1)}),(1_b^{(1)},2^{(1)},5^{(1)}), $ $
(1_c^{(1)},3^{(1)},4^{(1)}),(2^{(1)},5^{(1)},4^{(1)})\}$.

The second component, $H_2$, contains the set of vertices
$\{2^{(2)},3^{(2)},4^{(2)},5^{(2)}, 1_a^{(2)}, 1_c^{(2)}\}$, with the
vertex $1_a^{(2)}$ frozen to have spin $\sigma_1$ and the vertex
$1_c^{(2)}$ frozen to have spin $\sigma_0$; and the hyperedges
being $\{(1_a^{(2)},2^{(2)},3^{(2)}),$ $ (1,2^{(2)},5^{(2)}), $ $
(1_c^{(2)},3^{(2)},4^{(2)}), (2^{(2)},5^{(2)},4^{(2)})\}$.
Finally, the third component, $H_3$, contains the set of vertices
$\{2^{(3)},3^{(3)},4^{(3)},5^{(3)},  1_a^{(3)}, 1_b^{(3)}\}$; the
vertices $1_a^{(3)}$ and $1_b^{(3)}$ frozen to have spin
$\sigma_1$; and hyperedges $\{(1_a^{(3)},2^{(3)},3^{(3)}), $ $
(1_b^{(3)},2^{(3)},5^{(3)}), $ $
(1,3^{(3)},4^{(3)}),(2^{(3)},5^{(3)},4^{(3)})\}.$

It is clear that the following holds,
\begin{align*}
\frac{P_H(x_1 = \sigma_1)}{P_H(x_1= \sigma_0)} = ~ &
\frac{P_{H_1}(x_1 = \sigma_1)}{P_{H_1}(x_1= \sigma_0)} \times
\frac{P_{H_2}(x_1 = \sigma_1)}{P_{H_2}(x_1= \sigma_0)} \times
\frac{P_{H_3}(x_1 = \sigma_1)}{P_{H_3}(x_1= \sigma_0)} \\
= ~ & \frac{P_{G_1}(x_{1_a} = \sigma_1|
x_{1_b}=\sigma_0,x_{1_c}=\sigma_0)}{P_{G_1}(x_{1_a} = \sigma_0|
x_{1_b}=\sigma_0,x_{1_c}=\sigma_0)} \times \frac{P_{G_1}(x_{1_b} =
\sigma_1| x_{1_a}=\sigma_1,x_{1_c}=\sigma_0)}{P_{G_1}(x_{1_b} =
\sigma_0|
x_{1_a}=\sigma_1,x_{1_c}=\sigma_0)}  \\
& \quad \times \frac{P_{G_1}(x_{1_c} = \sigma_1|
x_{1_a}=\sigma_1,x_{1_b}=\sigma_1)}{P_{G_1}(x_{1_c} = \sigma_0|
x_{1_a}=\sigma_1,x_{1_b}=\sigma_1)}, \\
= ~ & \frac{P_G(x_1 = \sigma_1)}{P_G(x_1= \sigma_0)} .
\end{align*}

Further, this general procedure for separating the children of the
root can now be performed iteratively on each of its children to
yield a CD hypertree, $H_{CD}$, in the same way as one generates
the CD tree for pairwise interactions. Since at each stage, the
number of unfrozen vertices reduces by one, the procedure
terminates yielding a hypertree with the degree of every vertex
bounded by its degree in the original hypergraph.  This leads to
the following result for the case of hypergraphs,

\begin{theorem}
\label{th:comptreepresmargprobhypergraphs} For every hypergraph
$G=(V,E)$, every $\Lambda \subseteq V$, any configuration $
\sigma_\Lambda$, and all $\sigma_v$
\[ R_G^{\sigma_\Lambda}(v=\sigma_v) =
\mathbb{R}_{H_{CD}}^{\sigma_\Lambda}(v=\sigma_v), \] where
$\mathbb{R}_{H_{CD}}^{\sigma_\Lambda}(v=\sigma_v)$ stands for the
ratio (with respect to the reference spin, say $\sigma_0$) of the
probability that the root $V$ of $H_{CD}$ has spin $\sigma_v$ when
computations are performed as described previouly. The actual
probabilities can be computed from the ratios by normalizing them
such that the probabilities sum to one.
\end{theorem}

\section{Spatial mixing and Infinite regular trees}

\label{sse:ssm}

In this section, we study spatial mixing and demonstrate
sufficient conditions for spatial mixing to exist for all graphs
$G$ with maximum degree $b+1$ in terms of spatial mixing
conditions on the infinite regular tree, $\hat{\mathbb{T}}^b$, of
degree $b + 1$. We review the concept of strong spatial mixing
that was considered in \cite{wei06} and prove one of our main
results.

\begin{definition}
\label{def:ssm} Let $\delta: \mathbb{N} \to \mathbb{R}^+$ be a
function that decays to zero as $n$ tends to infinity. The
distribution over the spin system depicted by $G=(V,E)$ exhibits
{\em strong spatial mixing} with rate $\delta(\cdot)$ if and only
if for every spin $\sigma_1$, every vertex $v \in V$ and $\Lambda
\subseteq V$ and any two spin configurations, $\sigma_\Lambda,
\tau_\Lambda,$ on the frozen spins, we have
$$ |p(v=\sigma_1| X_\Lambda = \sigma_\Lambda) - p(v=\sigma_1| X_\Lambda =
\tau_\Lambda) | \leq \delta(\mathrm{dist}(v,\Delta)), $$ where
$\Delta \subseteq \Lambda$ stands for the subset in which the
frozen spins differ.
\end{definition}

Let $T$ denote a rooted tree. We say that a collection $L$ of virtual
edges is a set of  {\em valid coupling lines}, if they satisfy the
following constraints: a coupling line joins a vertex to some vertex
in the subtree under it; the lower endpoints of the coupling lines
are unique; no pair of coupling lines form a nested pair or an
interleaved pair, i.e. the endpoints do not lie on a single path.

\begin{rem}
Observe that the pruned CD tree, $T_{CD}^\Lambda$, is a tree with
a set of valid coupling lines. The pruned CD tree also has the
property that the end points of coupling lines have a
corresponding twin leaf that is frozen to $\sigma_q$, but we have
not imposed that requirement above. It is possible that enforcing
that requirement and thus limiting the set of {\em valid coupling
lines} may strengthen the results, but we omit it here for ease of
exposition.
\end{rem}

\begin{definition}
\label{def:vssm} Let $T$ denote a rooted tree. Let $\delta:
\mathbb{N} \to \mathbb{R}^+$ be a function that decays to zero as
$n$ tends to infinity. The distribution over the spin system at
the root, $v$, of $T$ exhibits {\em very strong spatial mixing}
with rate $\delta(\cdot)$ if and only if for every spin
$\sigma_1$, every set of {\em valid coupling lines},  for every
$\Lambda \subseteq V$ and any two spin configurations,
$\sigma_\Lambda, \tau_\Lambda,$ on the frozen spins, we have
$$ \Big| p_T(v=\sigma_1| X_\Lambda = \sigma_\Lambda) - p_T(v=\sigma_1| X_\Lambda =
\tau_\Lambda) \Big| \leq \delta(\mathrm{dist}(v,\Delta)), $$ where
$\Delta \subseteq \Lambda$ stands for the subset in which the frozen
spins differ. The computations of the marginal probability on this
tree with coupling lines is performed as described in Section
\ref{sse:compsec}.
\end{definition}

\begin{rem}
From the recursions observe that the computation tree can be
pruned at any frozen vertex or at any lower endpoint of a coupling
line.
\end{rem}

\begin{rem}
It is clear that very strong spatial mixing reduces to strong
spatial mixing in the absence of coupling lines. Thus very strong
spatial mixing on a tree implies strong spatial mixing with the
same rate on the tree.
\end{rem}

The main result of this section is that very strong spatial mixing
on the infinite regular tree with degree $b+1$ implies  strong
spatial mixing on any graph with degree $b+1$. We will distinguish
between two cases of neighboring interactions:
\begin{itemize}
\item[$(i)$]Spatially invariant interactions $\Phi(\cdot,\cdot) \geq
0$ and potentials $\phi(\cdot)\geq 0$ where the interaction matrix
$\Phi(\cdot,\cdot)$ satisfies the positively alignable condition
stated below.
\item[$(ii)$] General spatially invariant interactions $\Phi(\cdot,\cdot) \geq 0$ and
potentials $\phi (\cdot)\geq 0$ that need not satisfy the positively
alignable condition.
\end{itemize}

\begin{definition}
\label{def:posalign} A matrix $\Phi(\cdot,\cdot)$ is said to be
{\em positively alignable} if there exists a non-negative vector
$\alpha(\cdot)$ such that the column vectors of the matrix $\Phi$ can be
aligned in the $[1 ... 1]^T$ direction, i.e. $\Phi \alpha = [ 1 1 ...
1]^T$. Alternately, the vector $[1 ... 1]^T$ belongs to the convex
cone of the column vectors of $\Phi$.
\end{definition}

 Note that a sufficient
condition for $\Phi$ to be positively alignable is the existence
of a ({\em permissive} spin) $\sigma_0$ which satisfies the following property:
$\Phi(\sigma_i,\sigma_0) = c_1>0$ for all spins $\sigma_i$, and
$\phi(\sigma_0) = c_2>0$ (e.g. the ``unoccupied" spin in
independent sets).

\begin{rem}
We will state the next theorem for Case $(i)$, and a similar
theorem (see Section \ref{sse: geninter}) will hold for the other
case. The reason for separating the two cases is that in Case
$(i)$ one can stay within the same spin space in the infinite tree
$\hat{\mathbb{T}}^b$, to verify very strong spatial mixing.
\end{rem}

\subsection{Interactions that are positively alignable}

\begin{theorem}
\label{thm:treesuffices} For every positive integer $b$ and fixed
$\Phi(\cdot, \cdot),\phi(\cdot)$ such that $\Phi$ is positively
alignable, if $\hat{\mathbb{T}}^b$ exhibits very strong spatial
mixing with rate $\delta$ then every graph with maximum degree $b+1$
and having the same $\Phi(\cdot, \cdot),\phi(\cdot)$ exhibits strong
spatial mixing with rate $\delta$.
\end{theorem}

\begin{proof}
The proof of this theorem follows in a straightforward manner from
Theorem \ref{th:comptreepresmargprob}. If $T_\Lambda$ is the tree
in Section \ref{sse:comptree} rooted at $v$, i.e. $T_{CD}$ adapted
to $\Lambda$, then Theorem \ref{th:comptreepresmargprob} implies
that
\begin{equation}
\label{eq:gratree} \Big|p_G(v=\sigma_1| X_\Lambda =
\sigma_\Lambda) - p_G(v=\sigma_1| X_\Lambda = \tau_\Lambda) \Big|
= \Big|p_{T_\Lambda}(v=\sigma_1| X_\Lambda = \sigma_\Lambda) -
p_{T_\Lambda}(v=\sigma_1| X_\Lambda = \tau_\Lambda)\Big|.
\end{equation}
Further note that for any subset $\Delta$ of vertices of $G$,
dist(v,$\Delta$) is equal to the distance between the root $v$ and
the subset of vertices of $T_\Lambda$ composed of the copies of
vertices in $\Delta$ as the paths in $T_\Lambda$ correspond to paths
in $G$.  To complete the proof we need to move from $T_\Lambda$ to
$\hat{\mathbb{T}}^b$.

Note that $\Phi$ is positively alignable is equivalent to the
existence of a probability vector $a(\cdot) $ such that
\begin{equation}
\label{eq:posalign} \sum_i \Phi(\sigma_l,\sigma_i) \phi(\sigma_i)
a(\sigma_i) = c_1
>0, ~ \forall \sigma_l.
\end{equation}
As every vertex in $T_\Lambda$ has at most the degree of the
vertex in $G$, one can view $T_\Lambda$ as a subgraph of
$\hat{\mathbb{T}}^b$. (As before $\Lambda$ is also assumed to
contain the vertices that are frozen to $\sigma_q$ by the
construction.) Let $\partial(T_\Lambda)$ represent the non-fixed
boundary vertices, i.e. vertices in $T_\Lambda$ that are not fixed
by $\Lambda$, are not the lower end points of a dotted line, and
have degree strictly less than $b + 1$. Let $\Lambda_1$ denote the
set of vertices:\ in $\hat{\mathbb{T}}^b \setminus T_\Lambda$ that
is attached to one of the vertices in $\partial(T_\Lambda)$.
Append $\Lambda_1$ to $T_\Lambda$ to yield a subtree,
$\hat{\mathbb{T}}^b_\Lambda$  of $\hat{\mathbb{T}}^b$. Choose the
spins for the vertices in $\Lambda_1$ independently, distributed
proportional to $\phi(\cdot)a(\cdot)$.

We claim that
\[ p_{T_\Lambda}(v=\sigma_1| X_\Lambda = \sigma_\Lambda) =
p_{\hat{\mathbb{T}}^b_\Lambda}(v=\sigma_1| X_\Lambda =
\sigma_\Lambda). \] This follows from the observation that for all
$u_i$ in $\Lambda_1$ we have
\begin{equation*}
\begin{split}
& \frac{\sum_{l=1}^q \Phi_{v,u_i}(\sigma_v,\sigma_l)
R_{T_i}^{\sigma_{\Lambda_i}}(u_i=\sigma_l)}{\sum_{l=1}^q
\Phi_{v,u_i}(\sigma_q,\sigma_l)
R_{T_i}^{\sigma_{\Lambda_i}}(u_i=\sigma_l)} \\
& \quad = \frac{\sum_{l=1}^q \Phi_{v,u_i}(\sigma_v,\sigma_l)
a(\sigma_l)\phi(\sigma_l)}{\sum_{l=1}^q
\Phi_{v,u_i}(\sigma_q,\sigma_l) a(\sigma_l)\phi(\sigma_l)}
\stackrel{(a)}{=}  1,
\end{split}
\end{equation*}
where $(a)$ follows from \eqref{eq:posalign}. Thus the recursions
in $\hat{\mathbb{T}}^b_\Lambda$ becomes identical to the ones in
$T_\Lambda$.

Now from the very strong spatial mixing property that
$\hat{\mathbb{T}}^b$ is assumed to possess, we have
\begin{equation*}
\begin{split}
& \Big| p_{T_\Lambda}(v=\sigma_1| X_\Lambda = \sigma_\Lambda) -
p_{T_\Lambda}(v=\sigma_1| X_\Lambda = \tau_\Lambda) \Big| \\
& \quad= \Big| p_{\hat{\mathbb{T}}^b_\Lambda}(v=\sigma_1|
X_\Lambda = \sigma_\Lambda) -
p_{\hat{\mathbb{T}}^b_\Lambda}(v=\sigma_1| X_\Lambda =
\tau_\Lambda) \Big| \leq \delta(\mathrm{dist}(v,\Delta)).
\end{split}
\end{equation*}
The above equation along with \eqref{eq:gratree} completes the
proof.
\end{proof}

\begin{coro}
\label{cor:uniqueGM} Very strong spatial mixing on
$\hat{\mathbb{T}}^b$ (with positively alignable $\Phi$) implies a
unique Gibbs measure on all graphs with maximum degree $b+1$.
\end{coro}
\begin{proof}
From Theorem \ref{thm:treesuffices}, very strong spatial mixing on
$\hat{\mathbb{T}}^b$ with positively alignable $\Phi$ implies strong
spatial mixing on graphs with maximum degree $b+1$. Since  strong
spatial mixing is a sufficient condition for the existence of a
unique Gibbs measure on all graphs with maximum degree $b+1$, the
result follows.
\end{proof}

\subsection{General Interactions} \label{sse: geninter}

Consider the scenario of general interactions. Define an extra (permissive)
spin $\sigma_0$ that satisfies the following property:
$\Phi(\sigma_0,\sigma_l) = c_2 > 0, \phi(\sigma_0) = c_3 > 0.$ If
there is very strong spatial mixing on the infinite tree with this
extra spin $\sigma_0$ then the following analogue of Theorem
\ref{thm:treesuffices} holds.

\begin{theorem}
\label{thm:gentreesuffices} For every positive integer $b$,
 if $\hat{\mathbb{T}}^b$ (with the extra spin $\sigma_0$) exhibits very strong spatial
mixing with rate $\delta$ then every graph with maximum degree $b+1$
and having the same $\Phi(\cdot,\cdot),\phi(\cdot)$  exhibits very
spatial mixing with rate $\delta$.
\end{theorem}

\begin{proof}
The proof is similar to that of Theorem \ref{thm:treesuffices}
except for the following changes. Fix the spins of the vertices in
$\Lambda_1$ to $\sigma_0$ instead of generating them independently
with probability $a(\cdot)$. Condition also on the event that none
of the sites in $T_\Lambda$ are assigned the extra spin $\sigma_0$.
With these two changes made, the proof of Theorem
\ref{thm:treesuffices} carries over and hence is not repeated.
\end{proof}

\subsection{ On very strong spatial mixing on
trees} The idea of very strong spatial mixing is different from
the standard notions of spatial mixing due to the introduction of
coupling lines. However, it is key to note that these coupling
lines behave similarly in configurations $\sigma_\Lambda$ and
$\tau_\Lambda$ and thus conceptually it is similar to strong
spatial mixing where vertices close to the root are allowed to be
frozen to identical spins in both $\sigma_\Lambda$ and
$\tau_\Lambda$. However the fact that the actual computations
involve spins to be frozen to different values may lead to a
strictly stronger condition than strong spatial mixing. In some
sense, this condition demands that the difference of marginal
probabilities depend only on the spatial locations of the
frozen vertices and not on the spins that these vertices assume,
reminiscent of uniform convergence in analysis.

One sufficient condition for very strong spatial mixing is the
existence of a Lipschitz contraction for probabilities or
log-likelihoods, as in \cite{ban06,bgknt06}. In general if one can
show that some continuous monotone function
$f(p^{\sigma_\Lambda}(\sigma_v))$, where
$p^\sigma_\Lambda(\sigma_v)$ is computed using the recursions in
\eqref{eq:treerec0} from the probabilities of its children $\{
p_i^{\sigma_\Lambda}(\sigma_l) \}$, satisfies
\[ | f(p^{\sigma_\Lambda}(\sigma_v)) -
f(p^{\tau_\Lambda}(\sigma_v)) | < K \max_{i,l} |
f(p_i^{\sigma_\Lambda}(\sigma_l)) - f(p_i^{\tau_\Lambda}(\sigma_l))
|
\]
for some $K < 1$, then one can show that this implies very strong
spatial mixing (indeed with an exponential rate).

\section{Algorithmic implications}

The idea of strong spatial mixing, combined with an exponential
decay of correlation, has been used recently in
\cite{wei06,gak07,bgknt06} to derive polynomial time approximation
algorithms for counting problems like independent sets, list
colorings and matchings. Traditionally these counting problems
were approximated using Markov chain Monte Carlo  (MCMC) methods yielding
randomized approximation algorithms. In contrast the new
techniques based on spatial correlation decay yield {\em
deterministic} approximation algorithms, thus providing  a new alternative
to MCMC techniques.

\begin{definition}
\label{def:expdecay} A pairwise interacting system
($\Phi(\cdot,\cdot),\phi(\cdot)$) is said to have an {\em
exponential strong spatial correlation decay} if an infinite
regular tree of degree $D$, rooted at $v$, has a very strong
spatial mixing rate, $\delta(\mathrm{dist}(v,\Delta)) \leq
e^{-\kappa_D \mathrm{dist}(v,\Delta)}$ for some $\kappa_D > 0$.
\end{definition}

From the previous two sections, we will see that  the marginal
probabilities (and thus the partition function) for any pairwise
interacting system with finite spins with an exponential strong
spatial correlation decay, whose interactions can be modeled as a
graph $G$ with bounded degree, can be approximated efficiently.

\begin{lemma}
\label{lem:ptasexpdecay} Consider a graph $G$ of bounded degree,
say $D$, denoting the interactions of a pairwise interaction
system with  exponential strong spatial correlation decay. Then
the marginal probability of any vertex $v$ can be approximated to
within a factor $(1\pm \epsilon)$, for $\epsilon = n^{-\beta}$, in
a polynomial time given by $\Theta(n^{\frac{\beta}{\kappa_D}\log
D})$.
\end{lemma}

\begin{proof}
From the definition of strong spatial mixing rate it is clear that
the marginal probability at the root can be approximated to a
$(1+\epsilon)$ factor, provided $\mathrm{dist}(v,\Delta) > -
\frac{\log \epsilon} {\kappa_D} =: \ell$. That is, for any initial
assignment of marginal probabilities to lead nodes at depth $l$
from the root, the recursions in \eqref{eq:treerec} would give a
$(1+\epsilon)$ approximation to the true marginal probability.

Let $C$ denote the computation time required for one step of the
recursion in \eqref{eq:treerec}, then it is clear that  computing
the probability at the root given the marginal probabilities at
depth $\ell$ requires $\Theta([(q-1)D]^\ell)$ time. The hidden
constants in $\Theta$ depend on $C$ and $q$. Observe that a bound
for the computation time, $t_\ell$, at depth $\ell$ can be obtained
via the recursion $t_\ell \leq qC + (q-1) D t_{\ell - 1}$.

Therefore, if one wishes to obtain an $\epsilon = n^{-\beta}$
approximation, then the computational complexity would be
$\Theta(n^{\frac{\beta}{\kappa_D}\log (q-1)D}).$ Thus, the marginal
probability as well as the partition function can be approximated in
polynomial time.
\end{proof}

\begin{rem}
It is well known that the partition function can be computed as a
telescopic product of marginal probabilities (of smaller and
smaller systems) and thus an efficient procedure for yielding the
marginal probabilities also yields an efficient procedure (usually
time gets multiplied by $n$ and the error gets magnified by $n$)
for computing the partition function.
\end{rem}

\section{Remarks and conclusion}

\noindent {\em On colorings in graphs:} Consider
the anti-ferromagnetic hard-core Potts model with $q$ spins, or
equivalently, consider the vertex coloring of graph $G$ with $q$
colors. It is conjectured that for any infinite graph with maximum degree $D$
(and with appropriate vertex transitivity assumptions, so
that the notion of Gibbs measures make sense),
one can show that this system has a unique Gibbs
measure as long as $q$ is at least $D + 1$. Using the results in
the previous sections, if one establishes that the infinite
regular tree with degree $D$ has very strong spatial mixing when
$q$ is at least $D + 1$, then this will imply that any graph with
maximum degree $D$ will also have very strong spatial mixing
and thus a unique Gibbs measure.

It is known from \cite{jon02} that the infinite regular tree with
degree $D$ has {\em weak} spatial mixing when the number of colors is
at least $D+1$. The nature of the correlation decay suggests
that very strong spatial mixing should also hold in this instance.
However, it is not clear to the authors that the proof can be
modified to provide an argument for very strong spatial mixing (or
even whether the proof can be extended to show weak spatial mixing
for irregular trees with maximum degree $D$).
 Another possible
approach that is yet to be explored completely is whether weak
spatial mixing and some monotonicity arguments, like in
\cite{wei06} for independent sets, will directly imply very strong
spatial mixing.

\medskip

\noindent {\em Conclusion:} We have shown the existence of a
computation tree in graphical models that compute the exact
marginal probabilities in any graph. Further we have shown that
from the point of view of very strong spatial mixing, a notion of
spatial correlation decay, the infinite regular tree is a
worst-case graph. So proving results on infinite regular trees
would immediately imply similar results for graphs with bounded
degree.

\section*{Acknowledgements}
The authors would like to thank Mohsen Bayati, Christian Borgs,
 Jennifer Chayes,  Marc Mezard, Andrea
Montanari for helpful comments and useful discussions. Special
thanks go to  Elchanan Mossel for urging several of us, interested
in this topic, to work on this problem as well as for useful
discussions with the authors. The authors would like to thank Dror
Weitz for making useful comments and suggestions and for identifying
an error in an earlier version which has led us to redefine the very
strong spatial mixing condition on trees.

\bibliographystyle{amsalpha}

\bibliography{mybiblio}

\end{document}